# Quantum groups and $q$-lattices in phase space[*]

**J. Wess**

*Max-Planck-Institut für Physik*
*(Werner-Heisenberg-Institut)*
*Föhringer Ring 6, D-80805 München*

ABSTRACT

Quantum groups lead to an algebraic structure that can be realized on quantum spaces. These are noncommutative spaces that inherit a well defined mathematical structure from the quantum group symmetry. In turn such quantum spaces can be interpreted as noncommutative configuration spaces for physical systems which carry a symmetry like structure. These configuration spaces will be generalized to noncommutative phase space. The definition of the noncommutative phase space will be based on a differential calculus on the configuration space which is compatible with the symmetry. In addition a conjugation operation will be defined which will allow us to define the phase space variables in terms of algebraically selfadjoint operators. An interesting property of the phase space observables will be that they will have a discrete spectrum. These noncommutative phase space puts physics on a lattice structure.

---

Let me start by presenting a simple example of a $q$-deformed Lie algebra. This will be the algebra $sU_q(2)$ which is interesting by itself and exhibits many of the features that we encounter with quantum groups. I start with the defining algebra which for $q = 1$ is, up to a normalization, easily identified with the algebra $sU(2)$.

$$\frac{1}{q}T_+T_- - qT_-T_+ = T_3$$
$$q^2 T_3 T_+ - \frac{1}{q^2}T_+ T_3 = (q + \frac{1}{q})T_-$$
$$\frac{1}{q^2}T_3 T_- - q^2 T_- T_3 = -(q + \frac{1}{q})T_-$$

In general, $q$ is a complex parameter ($q \neq 0$). In this lecture I will assume $q$ to be real and, moreover $q \geq 1$. The algebra $sU(2)$ ($q = 1$) is treated in all elementary text books on quantum-mechanics as the algebra of angular momentum. Thus let me try to follow such a text book as far as possible for the case $q > 1$.

To assure real eigenvalues for angular momentum we demand conjugation properties, especially $\bar{T}_3 = T_3$.

The $q$-deformed algebra is consistent with the conjugation property:

$$\bar{T}_3 = T_3, \bar{T}_+ = \frac{1}{q^2}T_-$$

There is an operator that represents the square of the total angular momentum. It commutes with all the components of the angular momentum. Such an operator exists in the $q$-deformed case as well. to write it in a compact form we first introduce the notation:

$$\lambda = q - \frac{1}{q}, \tau = 1 - \lambda T_3$$

We find that

$$\mathcal{C} = q^2(T_- T_+ + \frac{1}{\lambda^2})\tau^{-\frac{1}{2}} + \frac{1}{\lambda^2}\tau^{\frac{1}{2}} - \frac{1+q^2}{\lambda^2}$$

commutes with $T_-, T_+$ and $T_3$ and that for $q = 1$, $\mathcal{C}|_{q=1} = \vec{j}^2$, if we denote angular momentum by $\vec{j}$. Physical states transform according to representations of the angular momentum.

The algebraic relations, conjugation and the Casimir operator $\vec{j}^2$ are all we need to construct such representations. With exactly the same methods, we construct the representations of $sU_q(2)$. They are 2j+1 dimensional and the states can be labelled by the eigenvalues of $T_3$. With the $q$-numbers

$$[n]_r = \frac{1 - q^{rn}}{1 - q^r}$$

the matrix elements are very similar to those of $sU(2)$:

$$T_3 |j, m\rangle = \frac{1}{q}[2m]_{-2}|j, m\rangle$$

$$T_+ |j, m\rangle = \frac{1}{q}\sqrt{[j+m+1]_{-2}[j-m]_2}|j, m+1\rangle$$
$$T_- |j, m\rangle = q\sqrt{[j+m]_{-2}[j-m+1]_2}|j, m-1\rangle$$
$$\mathcal{C}|j, m\rangle = [j]_{-2}[j+1]_2|j, m\rangle.$$

Except for a factor 2 that is related to the normalization of our generators, this representations tend to the familiar representations of angular momentum for $q \to 1$.

For convenience, we give the $j = \frac{1}{2}$ representation, the spinor representation of $sU_q(2)$:

$$T_+ |\frac{1}{2}, -\frac{1}{2}\rangle = \frac{1}{q}|\frac{1}{2}, \frac{1}{2}\rangle$$
$$T_+ |\frac{1}{2}, \frac{1}{2}\rangle = 0$$
$$T_- |\frac{1}{2}, -\frac{1}{2}\rangle = 0$$
$$T_- |\frac{1}{2}, \frac{1}{2}\rangle = q|\frac{1}{2}, -\frac{1}{2}\rangle$$
$$T_3 |\frac{1}{2}, -\frac{1}{2}\rangle = -q|\frac{1}{2}, -\frac{1}{2}\rangle$$
$$T_3 |\frac{1}{2}, \frac{1}{2}\rangle = \frac{1}{q}|\frac{1}{2}, \frac{1}{2}\rangle$$
$$\tau|\frac{1}{2}, \frac{1}{2}\rangle = q^2|\frac{1}{2}, \frac{1}{2}\rangle$$
$$\tau|\frac{1}{2}, -\frac{1}{2}\rangle = \frac{1}{q^2}|\frac{1}{2}, \frac{1}{2}\rangle$$

For $sU(2)$, starting with the spinor representation we can build all the other representations by "adding" angular momentum. Adding is based on the formula:

$$\vec{j} = \vec{j}_1 \otimes 1 + 1 \otimes \vec{j}_2$$

if $j_1$ and $j_2$ are the angular momenta that have to be added. The important fact is that $\vec{j}$ will satisfy the $sU(2)$ algebra if $\vec{j}_1$ and $\vec{j}_2$ do. The generalization of such a formula exists for $sU_q(2)$ as well. It is now called comultiplication

$$\Delta(T_3) = T_3 \otimes 1 + \tau \otimes T_3$$
$$\Delta(T_+) = T_+ \otimes 1 + \tau^{\frac{1}{2}} \otimes T_+$$
$$\Delta(T_-) = T_- \otimes 1 + \tau^{\frac{1}{2}} \otimes T_-$$

It is an easy exercise to verify that $\Delta(T)$ satisfy the algebra if the $T$'s do. The representations of the $sU_q(2)$ algebra as given above have the very interesting property that the non-vanishing matrix elements are at exactly the same entries of the matrix as those for $sU(2)$. As a consequence, the $T$ operators can be expressed in terms of the $j$ operators, or the $sU_q(2)$ operators are in the enveloping algebra of $sU(2)$.

To be more definite, if

$$j_+ j_- - j_- j_+ = j_0$$
$$j_0 j_\pm - j_\pm j_0 = \pm j_\pm$$
$$2j_+ j_- + j_0(j_0 - 1) = j(j+1)$$

we find that

$$T_3 = \frac{1}{\lambda}(1 - q^{-4j_0})$$
$$T_+ = \sqrt{2}q^{-j_0}\sqrt{\frac{[j_0+j]_{-2}[j_0-j-1]_2}{(j_0+j)(j_0+j-1)}}$$
$$T_- = q^2 \bar{T}_+$$

will represent the operators $T$.

Moreover, if $\bar{j}_+ = j_-$ and $\bar{j}_0 = j$ the $T$ operators will have the desired conjugation properties. Any operator that commutes with all the $j$ operators will commute with all the $T$ operators and vice versa.

We are now defining the concept of a quantum space. It is a space on which the $T$ operators act according to their representations and whose elements can be multiplied in accordance with the comultiplication rules of the $T$ operators.

For spinors:

$$T_+ x = qxT_+ + q^{-\frac{1}{2}}y$$
$$T_+ y = q^{-1}yT_+$$

$$T_- x = qxT_-$$
$$T_- y = q^{-1}yT_- + qx$$

$$T_3 x = q^2 xT_3 - qx$$
$$T_3 y = q^{-2} yT_3 + q^{-1}y$$

The inhomogeneous part on the right-hand side is copied from the spinor representation above and the homogeneous part reflects the comultiplication law.

As an example:

$$\Delta(T_3)|\frac{1}{2},\frac{1}{2}>\otimes|\rangle =$$
$$= -q|\frac{1}{2}\frac{1}{2}\rangle \otimes |\rangle + (\tau|\frac{1}{2}\frac{1}{2}\rangle) \otimes T_3|\rangle$$
$$= -q|\frac{1}{2}\frac{1}{2}>\otimes|\rangle + q^2|\frac{1}{2}\frac{1}{2}\rangle \otimes T_3|\rangle$$

If we dentote the first ket vector $|\frac{1}{2}\frac{1}{2}\rangle$ as $x$ and take away the second ket vectors $|\rangle$ we have

$$\Delta(T_3)x = -qx + q^2 xT_3.$$

All the actions of the $T$ operators on the $x\ y$ components are obtained the same way.

If we now have $T_+$ acting on $xy$ or $yx$ we find

$$T_+ xy = xyT_+ + q^{-\frac{1}{2}}y^2$$
$$T_+ yx = yxT_+ + q^{-\frac{3}{2}}y^2$$

It follows that

$$T_+(xy - qyx) = (xy - qyx)T_+$$

The relation

$$xy = qyx$$

is left invariant by $T_+$, and as can be easily seen, by all the other $T$ operators. The relation

$$xy = qyx$$

defines the two-dimensional quantum space of $sU_q(2)$. As harmless as this relation might look, it has a very interesting and non-trivial property, it does not generate any new relation in higher order. Instead of explaining this in more detail, I give another example that elucidates what I mean.

Try the relation

$$yx = xy + x^2 + y^2$$

and start from $y^2x$ which you want to reorder such that $x$ should appear at the left side. You will find

$$x^3 + y^3 + x^2y + xy^2 = 0$$

This is a third order relation. There are, as a consequence of the second order relation not as many linear independent polynomials of third degree as there are for commuting variables. This does not happen at any order if we start from the second order relation $xy = qyx$.

The interesting fact for us is that non-commutative spaces – a concept that we feel should be exploited for physics – arise in a natural way and with a well defined mathematical structure from quantum group symmetries.

Let us see how such a system can find a physical interpretation and what its consequences are. Starting from a non-commutative configuration space we will first generalize it to a phase space where non-commutativity is intrinsic for a quantum mechanical system. We thus shall regard the new non-commutative phase space as a $q$-deformation of the quantum mechanical phase space and apply all the ordinary machinery that we have learned from quantum mechanics for a physical interpretation. Thus we shall associate selfadjoint operators in a Hilbert space with observables and define the time-development of the system by a Schroedinger equation. These rules are not changed at all by $q$-deformation.

As an interesting result we shall see that the observables that in a natural way are associated with position and momentum, will have a discrete spectrum. Quantum groups put physics on a "$q$-lattice." We shall also see that this $q$-deformed phase space can be imbedded in an ordinary quantum mechanical phase space and that $q$-deformation has many similarities with a spontaneous breaking of a continuous spectrum to a discrete spectrum (continuous space to lattice).

The first step will be to base the $q$-deformed phase space on a non-commutative differential structure compatible with the non-commutative configuration space and the quantum group symmetry. In the next step we have to define conjugation properties on the non-commutative algebraic system that can be identified with hermitean conjugation of the corresponding operators in Hilbert space.

A consistent differential structure on the two component quantum space $(x, y)$ is:

$$\begin{aligned}\partial_x x &= 1 + q^2 x \partial_x + q\lambda y \partial_y \\ \partial_x y &= qy \partial_x \\ \partial_y x &= qx \partial_y \\ \partial_y y &= 1 + q^2 y\end{aligned}$$

We see that the ordinary Leibniz rule for differentiation had to be modified.

In this talk I would like to concentrate on the simplest possible $q$-deformed phase space structure. This would be a two-dimensional phase space where the Heisenberg commutation relations are $q$-deformed in accordance with our $q$-deformed differential structure ($-i\partial \to p$).

Thus I would like to ask you to accept as a model case for the deformation of the Heisenberg algebra the following algebra,

$$px - qxp = -1$$

Conjugation leads to

$$\bar{x}\bar{p} - q\bar{p}\bar{x} = i$$

It is clear that this does not allow both $x$ and $p$ to be mapped onto itself by conjugation. We could assume $\bar{x} = x$ and try to find $\bar{p}$ as a function of $p$ and $x$ such that the conjugate $q$-deformed Heisenberg relations are satisfied. To find this relation we first introduce a scaling operator $\Lambda$:

$$\Lambda = 1 + i(q-1)xp.$$

It has the property that:

$$\Lambda x = qx\Lambda \quad , \quad \Lambda p = q^{-1} p \Lambda$$

The occurrence of such an operator is typical for quantized spaces that emerge from quantum group symmetries. Note that for $q = 1$ follows $\Lambda = 1$, in the undeformed case, $\Lambda$ is not at our disposal.
It is now possible to verify that

$$\bar{p} = q^{-1}\Lambda^{-1} p$$

is a consistent conjugation assignment.
With this assignment follows

$$\bar{\Lambda} = q^{-1}\Lambda^{-1}.$$

To arrive at a hermitean momentum variable it is reasonable to define:

$$P = \frac{1}{2}(p + \bar{p})$$

To follow a convention we also define

$$X = \frac{1+q}{2q}x \qquad U = q^{\frac{1}{2}}\bar{\Lambda}$$

and we find Heisenberg relations for the $q$-deformed phase space

$$q^{\frac{1}{2}} XP - q^{-\frac{1}{2}} PX = iU$$
$$UX = q^{-1} XU, \, UP = qPU$$
$$\bar{P} = P, \bar{X} = X, \bar{U} = U^{-1}$$

For ordinary commutators it is necessary to introduce the factor $i$ on the right-hand-side of the commutator to be consistent with conjugation, for $q$-deformed commutators an operator has to appear. Here, and also for higher dimensions it is a unitary operator.

Hilbert space representations of this algebra are easily constructed. We first note that with $P, X$ also $\pi_0 P, \frac{1}{\pi_o} X$ will satisfy the algebra for arbitrary, real $\pi_0 \neq 0$.

We aim at representations where $P$ is represented by a selfadjoint linear operator and thus can be diagonalized. We assume $P$ to be diagonal and by a proper choice of $\pi_o$ we can always scale a non-vanishing eigenvalue to 1. From the algebra follows that $\overline{U}, U$ will change the eigenvalue by a factor $q, q^{-1}$, respectively. Thus we expect eigenvalues $q^n$, $n$ integer. A representation in which $\pi_0$ is adjusted to these eigenvalues shall be denoted by $\hat{P}_+, \hat{X}_+, \hat{U}_+$. Plus indicates that the eigenvalues of $\hat{P}_+$ are all positive. It is easy to see that the operators

$$\begin{aligned}\hat{P}_+ |n\rangle &= q^n \,|\, n\rangle, \quad \hat{U}_+ |n\rangle =|\, n-1\rangle \\ \hat{X}_+ |n\rangle &= i\frac{q^{-n}}{\lambda}\left\{q^{\frac{1}{2}}|n-1> -q^{-\frac{1}{2}}|n+1\rangle\right\} \\ \langle n|m\rangle &= \delta_{n,m}\end{aligned}$$

will represent the algebra. With the assumption that $P$ is selfadjoint and has one eigenvalue equal one, this representation is unique. However, $\hat{X}_+$ is hermitean but not selfadjoint. That has as a consequence that $\hat{X}_+$ is not diagonalizable in this represention. Actually it has to any complex number exactly one normalizable eigenvector. This means that $\hat{X}$ has a one parameter family of selfadjoint extensions. But none of these selfadjoint extensions will satisfy the $q$-deformed phase space algebra. By adding two representations with eigenvalue $+1$ and $-1(\pi_0 = -1)$ for $P$, it is, however, possible to find a representation where $P$ as well as $X$ are represented by selfadjoint linear operators. In this representation

$$\hat{X} = \hat{X}_+ \oplus \hat{X}_-$$

is diagonalizable by a $q$-deformed Fourier transformation which uses $q$-deformed cos and sin functions that have been defined by T.H. Kornwinder. The eigenvalues of $\hat{P}$ and $\hat{X}$ are both $\pm q^n$, $n$ integer. The corresponding eigenvectors will be denoted by $|n, \sigma>$, where $\sigma = +, -$, depending on the sign of the eigenvalue.

From these representations we can build new ones by adjoining a central element $\Pi$ to the algebra ($\Pi$ commutes with all other elements); with $P, X$, also $\Pi P, \pi^{-1} X$ will represent the algebra. We choose a representation for $\Pi$ with eigenvalues $s$, $1 \leq s < q$ (now we assume $q > 1$) and we normalise the eigenvectors by $\langle s' \mid s \rangle = \delta(s - s')$. The vectors $|n, \sigma\rangle \mid s_0 >$ will then belong to a representation of $P$ with eigenvalues $s_0 q^n$ ($s_0$ fixed). In this representation $X$ will have eigenvalues $s_0^{-1} q^n$. This representation space we call $\mathcal{H}_{s_0}$. Now we have a good control about all the representations of the $q$-deformed phase space algebra with selfadjoint operators $P$ and $X$.

The elements of the $q$-deformed algebra can also be expressed in terms of the undeformed algebra, whose elements are denoted by $p$ and $x$; they satisfy $[x, p] = i$. This is also a general feature of quantum group algebras. In our case we find that the following expressions for $P, X$ and $U$ will satisfy the $q$-deformed algebra if $x$ and $p$ satisfy the usual canonical commutation relations:

$$P = \frac{[2z - \frac{1}{2}]}{2z - \frac{1}{2}} p, \quad X = x, \quad U = q^{2z}$$

where $z = -\frac{i}{4}(xp + px)$.
If we take $p$ and $x$ in the usual quantum mechanical representation, $P$ and $X$ will be highly reduzible.

Of course, $x$ and $p$ can be any pair of canonical conjugate variables. We can change the assignment by a canonical transformation. An interesting class of such transformation is:

$$\tilde{p} = f(z)p, \tilde{x} = xf^{-1}(z)$$

As a consequence, $\tilde{z} = z$. The conjugation property is valid for the $\tilde{p}, \tilde{x}$ variables if

$$\bar{f}(\bar{z}) = f(z + \frac{1}{2}).$$

This, for instance is true for the function

$$f^{-1}(z) = \frac{[2z - \frac{1}{2}]}{2z - \frac{1}{2}}$$

It changes the assignment $X, P$ to the undeformed variables as follows:

$$\begin{align} X &= \tilde{x} \\ P &= \frac{[2\tilde{z} - \frac{1}{2}]}{2\tilde{z} - \frac{1}{2}} \tilde{p} \\ U &= q^{2z} \end{align}$$

If we take $p$ and $x$ (or $\tilde{p}$ and $\tilde{x}$) in the usual quantum mechanical representation, $P$ and $X$ will be highly reduzible. The expression obtained for $p, x(\tilde{p}, \tilde{x})$ in terms of $\hat{X}, \hat{P}, \hat{U}$ satisfies the canonical commutation relations, but $p, (\tilde{x})$ is not a selfadjoint operator and therefore not an observable. The formal expression of $p, x$ in terms of $\hat{X}, \hat{P}, \hat{U}$, obtained by formally inverting these relations, satisfies the canonical commutation relations, but $p$ is not a selfadjoint operator, it could not be an observable.

If we choose $P, X$ and $x$ to be observables, any of the representations spanned by the vectors $|n\sigma >| s_0 >$ with fixed value $s_0$ would describe one physical system. Thus measuring one of these observables would determine the value of $s_0$ and then the respective representation space $\mathcal{H}_{s_o}$. If our Hamiltonian, describing the time development of the system, is an expression in $P$ and $X$, then the full time development takes place in this representation space $\mathcal{H}_{s_o}$. In this representation $P, X$ and $x$ have a discrete spectrum with eigenvalues $s_0 q^n, s_0^{-1} q^n$ respectively. This situation resembles very much the situation of spontaneous symmetry breaking. It, however, takes place in the phase space and not in the configuration space. The expectation value of the central operator $\Pi$ determines the lattice structure of the system.

Let us have a look at the simplest Hamiltonian of the above type

$$H = \frac{1}{2} P^2$$

We know that its spectrum is of the form:

$$E_n = \frac{1}{2} s_0^2 q^{2n} \qquad -\infty \leq n \leq \infty \quad .$$

There is no lowest eigenvalue. This Hamiltonian can be expressed in terms of the ordinary canonical variables, with $q = e^h$:

$$H = \frac{2}{\lambda^2} \ p \ \frac{q + \frac{1}{q} - 2\cos((xp + px)h)}{1 + (xp + px)^2} \ p.$$

Though $p$ is not a selfadjoint operator in $\mathcal{H}_{s_0}$ the Hamiltonian will be. Note that the denominator as well as the numerator is always positive and does not become zero.

A Hamiltonian of this type, when considered as an operator in the physical Hilbert space of ordinary quantum mechanics is diagonalizable. We would have to discover that there is the $P, X$ algebra, find its representations and identify $H$ with $\frac{1}{2} P^2$. By restricting the domain of $H$ to $\mathcal{H}_{s_0}$ we have a spontaneous reduction to a lattice structure for $x$ on which $x$ is an observable.

For curiosity, let us have a look at the classical system corresponding to the above Hamiltonian. The equation of motion can be solved by a canonical transformation.

We find as a result:

$$x(t) = 2\sqrt{E} \ t \frac{1}{\text{Sinh}(h)} \left[ \frac{\text{Cosh}(h) - \cos(4Eth)}{1 + 16 E^2 t^2} \right]^{\frac{1}{2}}$$

For $h \to 0$ we recover the classical motion $x = \sqrt{2E}\, t$. It is interesting to note that for $h \neq 0$ and large $t$, the motion becomes periodical, the frequency increases with the energy $E$.

As a conclusion I would say that quantum groups lead to dynamical systems with interesting features. There is a whole class of Hamiltonians which look complicated when formulated in terms of the usual quantum mechanical variables but become rather simple when formulated in terms of the $q$-deformed variables. This is not restricted to a one-dimensional phase space but works in three and more dimensions as well. For systems, where the classical motion can be solved or the Hamiltonian diagonalized due to a dynamical symmetry, an analogous system can be found in $q$-deformed phase space and solved by a $q$-deformed dynamical symmetry. Here we have discussed the simplest model of this type – free motion in one dimension. The harmonic oscillator in one or three dimensions or the hydrogen atom in three dimensions are other examples.

A general feature is that the interactions modify the phase space structure and lead to a $q$-lattice for appropriate phase space variables. In a way the structure of the space is determined by the interaction.